%% file: mobius_latex.tex
\documentclass[11pt]{amsart}

\issueinfo{00}{0}{Month}{Year}
\copyrightinfo{2008}{University of Calgary}
\pagespan{1}{2}
\PII{ISSN 1715-0868}
\include{cdmstdcmds}

\def\divides{\mid}

\def\bottom{\widehat 0}

\def\Int{\hbox{\rm Int}}

\def\ideal#1{\langle #1\rangle}
\def\F{{\mathcal F}}

\def\CC{\mathbb{C}}
\def\NN{\mathbb{N}}

\begin{document}

\title{An uncertainty principle for M\"obius inversion on posets}

\author{Marcel K. Goh}
\address{Department of Mathematics and Statistics, McGill University, Montr\'eal, Qu\'ebec, Canada}
\email{marcel.goh@mail.mcgill.ca}

%

\date{(date1), and in revised form (date2).}

\thanks{The author is funded by the Natural Sciences and Engineering Council of Canada.}

\begin{abstract}
We give conditions for a locally finite join-semilattice $P$ to have the property that
for any functions $f:P\to \CC$ and $g:P\to \CC$ not identically zero and
linked by the M\"obius inversion formula, the support of at
least one of $f$ and $g$ is infinite. This generalises
and gives an entirely poset-theoretic proof of a result of Pollack. Various examples and non-examples
are discussed.
\end{abstract}

\maketitle

\section{Introduction}

In harmonic analysis,
the celebrated uncertainty principle states that a function and its Fourier
transform cannot both have small support, where the notion of ``small'' is made precise in various ways
by different theorems.
A 2011 paper of P.~Pollack proves an analogue of this principle in which the Fourier
transform
is replaced by the number-theoretic M\"obius transform~\cite{pollack2011}.
Concretely, it is shown that the support of an arithmetic function
and that of its M\"obius transform cannot both be finite. The proof is short and uses basic properties of
complex power series. Stronger theorems related to the asymptotic densities of the functions involved
were subsequently proved by P.~Pollack and C.~Sanna~\cite{pollacksanna}

One is led to wonder whether Pollack's original result, which can be
expressed in entirely poset-theoretic terms, has a purely poset-theoretic proof.
In this note, we show that such a proof does indeed exist,
and that it is general enough to also apply to the poset of all subsets of $\NN$ ordered by set inclusion.

\section{The main theorem}

In this section we establish definitions and notation before proceeding with the proof of our main theorem.
Details concerning these definitions can be found in any introductory text on enumerative
combinatorics; see, e.g., R.~P.~Stanley's textbook~\cite{enumcomb}.

Our main object of interest is a poset $(P,\le)$. An {\it interval} in a poset is a subset of
the form $[x,y] = \{ z\in P : x\le z\le y\}$ for some elements $x,y\in P$.
A poset is said to be {\it locally finite}
if every interval in it is finite. An element $\bottom$ such that $\bottom \le x$ for all $x\in P$
is called a {\it bottom element}. For any element $x\in P$, the {\it principal order ideal generated by $x$}
is the set $\{y\in P : y\le x\}$. This is denoted by $\ideal x$ for short.

An element $z\in P$ is said to be an {\it upper bound} of elements $x$ and $y$ in $P$ if $z\ge x$ and $z\ge y$.
If such a $z$ has the property that for all upper bounds $z'$ of $x$ and $y$, $z\le z'$, then $z$ is uniquely
defined, denoted by $x\vee y$, and said to be the {\it join} of $x$ and $y$. This operation is commutative
and associative, so we may speak of joins of any finite number of elements. A poset is called a
{\it join-semilattice} if it contains a join of any finite set of its elements.

We say that an element $z\in P$ {\it covers} an element $y\in P$ if $z\ge y$ and $y \le z' < z$ implies that
$z' = y$. If the poset $P$ has a bottom element $\bottom$, then any element that covers $\bottom$ is called an
{\it atom}.

The {\it M\"obius function} $\mu_P$ of a poset $P$ with bottom element
is a function from the intervals of $P$ to the complex numbers obtained by setting
$\mu_P([x,x]) = 1$ for all $x\in P$ and recursively putting
$$\mu_P([x,y]) = -\sum_{x\le z<y} \mu(x,z)$$
for all $x\le y$ in $P$. For brevity we write $\mu_P(x,y)$ instead of $\mu_P([x,y])$. The M\"obius inversion theorem
states that if $f$ is a function from $P$ to $\CC$ and $g:P\to \CC$ is given by
$$ g(y) = \sum_{x\le y} f(x)$$
for all $y\in P$, then for all $y\in P$ we have
$$ f(y) = \sum_{x\le y} \mu_P(x,y) g(x).$$
We shall say in this case that $g$ is the {\it M\"obius transform} of $f$.

Last but not least, for a function $f:P\to \CC$ the {\it support} of $f$ is defined to be the set
of $x\in P$ such that $f(x)\ne 0$.
Our main theorem gives conditions on a poset under
which $f$ and $g$ defined as above cannot both have finite support.

\begin{theorem}\label{main}
Let $P$ be a locally finite
join-semilattice with bottom element. Suppose that for every $y\in P$,
there exist infinitely
many $z\in P$ such that
\medskip
\item{i)} $z$ covers $y$; and
\smallskip
\item{ii)} for all $x\le y$, $\mu_P(x,z) = \mu_P(x,y)\mu_P(y,z) = -\mu_P(x,y)$.

\medskip
\noindent Then for any $f:P\to \CC$ that is not identically zero,
the support of $f$ and the support of its M\"obius transform $g$ cannot both be finite.
\end{theorem}

\begin{proof}
If $g$ has infinite support we are done, so suppose that $g$ has finite support
and choose $y$ with $f(y)\ne 0$. Let $S$ be the support of $g$ and let $T = S\cup\{y\}$. Let $m = \bigvee T$
be the join of these elements. Let $z$ be an element of $P$ that covers $y$ and note that if
$\ideal z\setminus \ideal y$ contains an element $x\in S$, then $z = x\vee y$, so $z\in \ideal m$.
We have shown that if $z\notin \ideal m$, then $S\cap \ideal z \setminus \ideal y = \emptyset$.

Since $P$ is locally finite with a bottom element, $\ideal m$ is finite; thus
by our hypothesis there exist infinitely many $z$ satisfying both conditions (i) and (ii)
that lie outside of $\ideal m$. Fix an arbitrary such $z$ and expand
\begin{align*}
f(z) &= \sum_{x \in \ideal z} \mu_P(x,z) g(x) \cr
&= \sum_{x\in \ideal y}\mu_P(x,z)g(x)
+ \sum_{x\in \ideal z \setminus \ideal y} \mu_P(x,z) g(x).\cr
\end{align*}
By the observation just proved, the second summation is zero, and by condition (ii),
$$f(z) = \sum_{x\in \ideal y} \mu_P(x,y)\mu_P(y,z) g(x) = \mu_P(y,z) f(y) = -f(y) \ne 0.$$
Hence we see that $f(z) \ne 0$ for infinitely many $z$.
\end{proof}

\section{Examples}

Let us say that a poset $P$ {\it adheres to the uncertainty principle} if for all pairs of functions
$f,g : P\to \CC$ linked by the M\"obius inversion formula, at most one of $f$ and $g$ has finite support.
In this section we show that our theorem applies not only to the divisibility poset treated by Pollack's
paper, but also to the poset of all finite subsets of $\NN$.

\medskip\noindent{\bf The divisibility poset.}\enspace
Let $P$ be the set $\NN$ of natural numbers,
partially ordered by divisibility.
The M\"obius function is $\mu_P(x,y) = \mu(y/x)$, where $\mu :\NN \to \{-1,0,1\}$ is given by
$$ \mu(n) = \begin{cases}
(-1)^k, & \textrm{if $n$ is the product of $k$ distinct primes;}\cr
0, & \textrm{if $n$ is divisible by a perfect square.}\cr
\end{cases}$$
This is a {\it multiplicative} function, in the sense that $\mu(mn) = \mu(m)\mu(n)$ whenever
$\gcd(m,n) = 1$. We use this in our proof of Pollack's result on arithmetic functions.

\begin{proposition}{\rm({\it Pollack,} 2011)}\label{primes}
Let $f:\NN\to \CC$ be a function not identically
zero and let $g:\NN\to \CC$ be given by
$$g(n) = \sum_{d\divides n} f(d),$$
so that
$$f(n) = \sum_{d\divides n} \mu(n/d) g(d).$$
Then the support of $f$ and the support of $g$ cannot both be finite.
\end{proposition}

\begin{proof}
Let $y\in \NN$ be given. There are infinitely many $z$ covering $y$, as we can let $z = qy$ for any prime
$q$. For any such $z$ and for all $d$ dividing $y$,
\begin{align*}
\mu_P(d,y)\mu_P(y,z) &= \mu(y/d)\mu(z/y) \cr
&= \mu(y/d)\mu(q) \cr
&= \mu(qy/d) \cr
&= \mu(z/d) \cr
&= \mu_P(d,z),
\end{align*}
since $\mu$ is a multiplicative function. Lastly, note that for each $z$ covering $y$,
$\mu(z/y) = \mu(q) = -1 \ne 0$. Thus by Theorem~\ref{main} we are done.
\end{proof}

Our proof only requires elementary number theory and the properties of posets, but a downside is that
it uses the infinitude of primes. Hence we cannot use Proposition~\ref{primes} to give
an alternative proof that there are infinitely many primes, as Pollack did in his paper.

\medskip\noindent{\bf Finite multisets of a countably infinite set.}\enspace
A {\it multiset on a ground set $X$} is a function $m : X\to \NN$. A multiset $m$ is said to be {\it finite}
if $\sum_{x\in X} m(x) <\infty$. We can place a partial order on the set of all finite multisets of $X$
by saying that $m\le m'$ if for all $x\in X$, $m(x)\le m'(x)$. When $X$ is countably infinite, this
is isomorphic to the divisibility poset above, since every integer $n$ can be uniquely represented by
the multiset $m_n$ mapping a prime $p$ to the exponent of $p$ in the factorisation of $n$. Thus we see
that the poset of all finite multisets on a countably infinite set adheres to the uncertainty principle.

\medskip\noindent{\bf Finite sets of natural numbers.}\enspace
Having dealt with multisets, we now turn to ordinary sets.
Let $P$ be the collection $\F$ of all finite subsets of $\NN$, ordered
by divisibility. In this context the M\"obius function $\mu_P$ is given by $\mu_P(S, T)  = (-1)^{|T|-|S|}$
for all $S\subseteq T$. It is not difficult to show that this poset satisfies the hypotheses of
Theorem~\ref{main}.

\begin{proposition}\label{subsets}
Let $f:\F\to \CC$ and let $g:\F\to \CC$ be given by
$$g(y) = \sum_{x\subseteq y} f(x),$$
so that
$$f(y) = \sum_{x\subseteq y} (-1)^{|y|-|x|} g(x).$$
Then the support of $f$ and the support of $g$ cannot both be finite.
\end{proposition}

\begin{proof}
Let $y\in \F$ be arbitrary. There are infinitely many $z$ covering $y$, since we can take
$y\cup \{n\}$ for any $n\in \NN\setminus y$, and for {\it any} chain of
sets $x\subseteq y\subseteq z$, we have
$$\mu_P(x,y)\mu_P(y,z) = \mu_\F(-1)^{|y|- |x|} (-1)^{|z|-|y|} = (-1)^{|z|-|x|} = \mu_P(x,z).$$
Hence $P$ adheres to the uncertainly principle, by Theorem~\ref{main}.
\end{proof}

Of course there is nothing special about $\NN$ here; our proof applies just as well to any other
countably infinite set.

\section{Non-examples}

It is worth noting that many infinite posets do not adhere to the uncertainty principle. Perhaps the simplest
example is given by the set $\NN$ of natural numbers under the usual definition of $\le$.
On this poset it is easy to check that
if we let $g(1) = 1$, $g(2) = -1$, and $g(n) = 0$ for all $n\ge 3$, then for $f(n) = \sum_{m\le n} g(n)$ we
have $f(1) = 1$ and $f(n) = 0$ for all $n\ge 2$.

The following proposition gives a simple necessary condition for a poset $P$ to adhere to the uncertainty
principle.

\begin{proposition}\label{necessary}
For a locally finite poset $P$ to adhere to
the uncertainty principle it is necessary that for all $x\in P$, the set
$$S_x = \{y \in P : \mu_P(x,y) \ne 0\}$$
is infinite.
\end{proposition}

\begin{proof}
Suppose there is $x\in P$ such that $S_x$ is finite. Let
$$g(y) = \begin{cases}
1, & \textrm{if $y=x$}; \cr 0, & \textrm{otherwise}.
\end{cases}$$
Clearly the support of $g$ is finite. Then since
$$f(y) = \sum_{x \le y} \mu_P(x,y) g(x) = \mu_P(x,y)$$
for all $y\ge x$ and $f(y) = 0$ for $y\not\ge x$, we conclude that the support of $f$ is contained in $S_x$.
Hence $P$ does not adhere to the uncertainty principle.
\end{proof}

Thus the fact that $P = (\NN,\le)$
does not adhere to the uncertainty principle is explained by the formula
$$\mu_P (m,n) = \begin{cases}
1, & \textrm{if $m=n$};\cr -1, & \textrm{if $m+1=n$};\cr 0, & \textrm{otherwise},
\end{cases}$$
which implies that, for example, $\mu_P(1,x)$ is nonzero only for $x\in \{1,2\}$.

\medskip\noindent{\bf Arbitrary convolutions.}\enspace
Let $\Int(P)$ denote the set of all intervals in a poset $P$. When $P$ is locally finite we define the
{\it incidence algebra} to be the algebra (over the complex field) of all functions from $\Int(P)$
to $\CC$, with multiplication given by the convolution operation
$$ (\alpha * \beta)(x,y) = (\alpha *\beta)\bigl([x,y]\bigr) = \sum_{x\le z\le y} \alpha(x,z)\beta(z,y).$$
The identity element of the incidence algebra is the function $\delta_P(x,y)$ defined by
$$\delta_P(x,y) = \begin{cases}
1, & \textrm{if $x = y$};\cr 0, & \textrm{otherwise.}
\end{cases}$$
The constant function $\zeta_P(x,y) = 1$ is called the {\it zeta function}, and the M\"obius inversion
formula is equivalent to the statement that $\mu_P$ and $\zeta_P$ are inverses to one another
in the incidence algebra; that is, $\mu * \zeta = \zeta * \mu = \delta$, where here (and in the rest
of the section) we dispense with
subscripts for brevity. (When $P$ has bottom element $\bottom$ we associate to any function $f:P\to\CC$
the function $f:\Int(P)\to \CC$ given by $f(\bottom, x) = f(x)$ and $f(z,x) = 0$ for all $z\ne 0$.)

A 2014 paper of C.~Sanna~\cite{sanna2014} generalises the density result of Pollack and Sanna found
in~\cite{pollacksanna} to arbitrary Dirichlet convolutions. In this vein, we formulate a more
general definition of our own. We shall say that $P$ adheres to the uncertainty
principle {\it with respect to $(\alpha,\beta)$} if $\alpha$ and $\beta$ are inverses of one another
in the incidence algebra of $P$ and whenever $f:P\to \CC$ and $g:P\to\CC$ are functions linked by
the relations
$$f(y) = \sum_{x\le y} \alpha(x,y)g(x)\qquad\hbox{and}\qquad
g(y) = \sum_{x\le y} \beta(x,y)f(x),$$
at most one of $f$ and $g$ can have finite support. Thus our previous definition of adhering
to the uncertainty principle is the same as doing so with respect to $(\mu,\zeta)$ under our new definition.

It is trivial to adapt the proof of Proposition~\ref{necessary} to the more general case.

\begin{proposition}\label{general}
Let $P$ be a locally finite poset and let $\alpha$ and $\beta$
be inverse to one another in the incidence algebra of $P$. Then for $P$ to adhere to the uncertainty
principle with respect to $(\alpha,\beta)$, it is necessary that for all $x\in P$, the sets
$$S_x = \{y \in P : \alpha(x,y) \ne 0\}\qquad\hbox{and}\qquad
T_x = \{y\in P : \beta(x,y)\ne 0\}$$
are both infinite.
\end{proposition}

\begin{proof}
Let $x\in P$ and suppose one of the sets above is finite. Define the function $f$ or
$g$ (depending on whether $T_x$ or $S_x$ is finite) accordingly and proceed as
in the proof of the previous proposition.
\end{proof}

Observe that
Proposition~\ref{necessary} only had to consider $\mu$ because $\zeta$ has infinite support for any infinite
poset.
The necessary condition put forth by Proposition~\ref{necessary}
is of course subsumed by the covering criterion
found in Theorem~\ref{main}, since any $z$ covering $y$ necessarily has $\mu_P(y,z) = -1\ne 0$.
One wonders whether the covering condition is necessary, and furthermore if the multiplicativity
condition (ii) in that theorem can be removed. We formulate the following hopeful conjecture in the broader
context pertaining to Proposition~\ref{general}.

\begin{guess}\label{conj}
Let $P$ be a locally finite join-semilattice with bottom element and let
$\alpha$ and $\beta$ be inverses to one another in the incidence algebra of $P$. Then $P$
adheres to the uncertainty principle with respect to $(\alpha,\beta)$ if and only if
for all $x\in P$, the sets
$$S_x = \{y \in P : \alpha(x,y) \ne 0\}\qquad\hbox{and}\qquad
T_x = \{y\in P : \beta(x,y)\ne 0\}$$
are both infinite.
\end{guess}
\goodbreak

\section*{Acknowledgements}

The author wishes to thank Carlo Sanna for suggesting the generalisation of Proposition~\ref{necessary},
which places Conjecture~\ref{conj} in a more general setting.

\end{document}

%% file: cdmstdcmds.tex
\usepackage{graphicx}
\usepackage{epstopdf}
\usepackage{amsthm,amsmath,amssymb,amsxtra,amscd}
\usepackage{verbatim}
\usepackage{rotating}
\usepackage{multirow}
\usepackage{multicol}
\usepackage{url}

\newtheorem{theorem}{Theorem}[section]

\newtheorem{guess}[theorem]{Conjecture}

\newtheorem{proposition}[theorem]{Proposition}

\newtheorem*{example*}{Example}

\newtheoremstyle{myexample}{3pt}{3pt}{\rmfamily}{}{\itshape}{:}{ }{\thmname{#1}\thmnumber{ #2}\thmnote{ (#3)}}
\theoremstyle{myexample}

\newtheoremstyle{myremark}{3pt}{3pt}{\rmfamily}{}{\itshape}{:}{ }{\thmname{#1}}
\theoremstyle{myremark}

\newtheorem*{observation*}{Observation}

\newtheoremstyle{conjecture}{3pt}{3pt}{\itshape}{}{\bfseries}{.}{ }{\thmname{#1}\thmnote{ (#3)}}
\theoremstyle{conjecture}

\newtheorem*{question*}{Question}

\newtheorem{theorem*}{Theorem}

\numberwithin{equation}{section}

\newcounter{algorithm}
\renewcommand{\thealgorithm}{\thesection.\arabic{algorithm}}
